\documentclass[12pt]{article}
\usepackage{graphics,latexsym,amssymb,amsmath,amscd}
\RequirePackage{graphics}
\usepackage[all,cmtip]{xy}
\usepackage[svgnames,table]{xcolor}
\def\abs#1{\left \vert #1 \right \vert} 
\def\[#1\]{\begin{eqnarray*}#1\end{eqnarray*}}

\def\SU{\hbox{\bf SU}}
\def\SL{\hbox{\bf SL}}
\def\phi{\varphi}

\input xy
\xyoption{all}
\newtheorem{thm}{Theorem}[section]
\newtheorem{dfn}[thm]{Definition}
\newtheorem{rmk}[thm]{Remark}
\newtheorem{cor}[thm]{Corollary}
\newtheorem{prop}[thm]{Proposition}
\newtheorem{lemma}[thm]{Lemma}

\newcommand{\Pf}{{\em Proof}. }
\newcommand{\EPf}{\hbox{}\hfill$\Box$\vspace{.5cm}}

\newcommand{\R}{{{\mathbb R}}}
\newcommand{\Z}{{{\mathbb Z}}}

\def\SU{\hbox{\bf SU}}
\def\SL{\hbox{\bf SL}}
\def\PSL{\hbox{\bf PSL}}

\def\phi{\varphi}

\def\R{{\mathbb R}}

\date{}
\title{Cartan connections and path structures with large automorphism groups}
\author{E. Falbel, M. Mion-Mouton and
 J. M. Veloso}
 
\newcommand{\Lm}{\mathcal{L}}
\newcommand{\Tm}{\mathcal{T}}
\newcommand{\Fitan}[1]{\ensuremath{\mathrm{T}#1}}
\newcommand{\cc}{\mathcal{C}}
\newcommand{\Diff}[2]{\mathopen{}\mathrm{D}_{#1}#2}
\newcommand{\norme}[1]{\ensuremath{\left\lVert#1\right\rVert}}
\newcommand{\restreinta}{\ensuremath{\mathclose{}|\mathopen{}}}
\newcommand{\Fiunitan}[1]{\ensuremath{{\mathrm{T}}^1{#1}}}
\newcommand{\Hn}[1]{\mathbf{H}^{#1}}
\usepackage{mathtools} 
\newcommand{\enstq}[2]{\ensuremath{\left\{#1\mathrel{}\middle|\mathrel{}#2\right\}}}
\DeclareMathOperator{\Aut}{Aut}
\newcommand{\Tan}[2]{\ensuremath{\mathrm{T}_{#1}#2}}
\DeclareMathOperator{\Lin}{L}
\DeclareMathOperator{\Ker}{Ker}
\newcommand{\ensemblenombre}[1]{\ensuremath{\mathbb{#1}}}
\newcommand{\N}{\ensemblenombre{N}}
\DeclareMathOperator{\id}{id}
\DeclareMathOperator{\Vect}{Vect}
\newcommand{\intervalle}[4]{\ensuremath{\mathopen{#1}#2
		\mathclose{}\mathpunct{};#3
		\mathclose{#4}}}
\newcommand{\intervalleff}[2]{\intervalle{[}{#1}{#2}{]}}
\usepackage[pdftoolbar=true,colorlinks=false]{hyperref}
\usepackage{todonotes}

\begin{document}
\maketitle
\newcommand{\D}{\mbox{$\cal D$}}
\newtheorem{df}{Definition}[section]
\newtheorem{te}{Theorem}[section]
\newtheorem{co}{Corollary}[section]
\newtheorem{po}{Proposition}[section]
\newtheorem{lem}{Lemma}[section]
\newtheorem{rem}{Remark}[section]
\newcommand{\Ad}{\mbox{Ad}}
\newcommand{\ad}{\mbox{ad}}
\newcommand{\im}[1]{\mbox{\rm im\,$#1$}}
\newcommand{\bm}[1]{\mbox{\boldmath $#1$}}
\newcommand{\sime}{\mbox{sim}}

\begin{abstract}
We classify compact manifolds of dimension  three equipped with a path structure and a fixed contact form 
(which we refer to as a strict path structure) 
{under} the hypothesis that their automorphism group is non-compact.  
{We use a Cartan connection associated to the structure and show that its curvature is constant.}
\end{abstract}

\section{Introduction}
\par The general context of this paper is the quest to understand, and eventually classify, 
manifolds equipped with both a geometric structure and a sufficiently large group of diffeomorphisms preserving that structure.   
One should mention here  Zimmer's program to classify actions of infinite groups which preserve volume, 
and the vague general conjecture appearing in \cite{DG}, stating that 
compact manifolds with both a rigid geometric structure and a sufficiently large group of automorphisms 
can be almost classified.

Several instances of this topic were studied, {in different geometric contexts and by numerous authors}. As an early  example we can mention Ferrand-Obata (\cite{ Fe, O}) theorem classifying conformal manifolds with non-compact groups of automorphisms  which was generalized to other rank one geometries in \cite{F0}. 
Another example is Ghy's classification of Anosov flows with smooth stable and unstable lines on compact three-manifolds (\cite{G}, see more details below).  This was generalised to higher dimensions  in \cite{BFL} using Gromov's theory of rigid geometric structures.  As a last example, among others, we mention the classification in \cite{Z} of compact Lorentz three-manifolds admitting a Killing flow not preserving any Riemannian metric.  

In the present work we used a description of the strict path geometric structure (see below) as a Cartan geometry.  The fact that one has a Cartan's connection which is invariant under an automorphism group with a dense orbit  implies that some components of its curvature vanish.  This will allow us to classify all such spaces. 

\subsection{Strict path structures with a non-compact automorphism group}
The particular instance of  the above problem we analyze is that of  three manifolds equipped with a so called 
Lagrangian contact structure or \emph{path structure}.  
This geometric structure (see \cite{IL} for a detailed introduction) 
has been studied since a long time (see, in particular, Cartan's study \cite{Ca}).
{It is related to the geometry of second order Ordinary Differential Equations (see \cite{BGH,IL}). {It is interesting to note that the first published paper by S. S. Chern (\cite{Ch}) concerns a generalization of this study to the case of third-order ODE.}}  
\par A \emph{path structure} on a three-manifold $M$ is a couple of two (smooth) one-dimensional distributions $(E^1, E^2)$ 
such that  $E^1\oplus E^2$ is a contact distribution.  See \cite{FV} for a theory corresponding to a complexification of this structure.
There are no obstructions to the existence of such a structure on a given compact orientable manifold, 
and the model space is the flag space of all sequences $\{V_1\subset V_2\subset \R^3\}$, 
where $V_1$ is a line and $V_2$ is a plane.  
This model space can be thought as the space of lines with a marked point in the projective plane $\R P^2$, 
or simply as the homogeneous space $\SL(3,\R)/B$ where $B$ is the Borel subgroup of upper triangular matrices.   

If, moreover, we fix a contact form $\theta$ whose kernel is the contact distribution $E^1\oplus E^2$, 
then we refer to the triplet $\Tm=(E^1,E^2,\theta)$ as a \emph{strict path structure}.
The flat model for strict path geometries is the Heisenberg space $\mathbf{Heis}(3)$ 
with two left-invariant directions and a fixed left-invariant contact form,
described in the paragraph \ref{subsection-flatstrict}.   
Its automorphism group is $\mathbf{Heis}(3)\rtimes P$, where $P$ is a group isomorphic to $\R^*$.   
On the other hand, 
a (non-flat) constant curvature model is given by a left-invariant structure on $\widetilde{\SL}(2,\R)$
(the universal cover of ${\SL}(2,\R)$, see section \ref{subsection-SL2strict}),
whose automorphism group is $\widetilde{\SL}(2,\R)\times\tilde{A}$, $\tilde{A}$ being a group isomorphic to $\R^*$.  
The main result of this paper is the following.
\begin{thm}\label{thm:GhyssansAnosov}
Let $\Tm$ be a strict path structure on a compact connected three-dimensional manifold $M$,
whose automorphism group is non-compact for the compact-open topology.
Then if $\Tm$ is of class $\cc^3$,
or if $\Tm$ is of class $\cc^2$ and has
a dense $\Aut^{loc}(M,\Tm)$-orbit:
\begin{enumerate}
\item either $(M,\Tm)$ is, up to a constant multiplication of its contact form,
isomorphic to $\Gamma\backslash\widetilde{\SL}(2,\R)$
for some discrete subgroup $\Gamma$ of $\widetilde{\SL}(2,\R)\times\tilde{A}$
acting freely, properly and cocompactly on $\widetilde{\SL}(2,\R)$;
\item or $(M,\Tm)$ is, up to a finite covering, isomorphic to $\Gamma\backslash\mathbf{Heis}(3)$ 
for some cocompact lattice $\Gamma$ of $\mathbf{Heis}(3)$.
\end{enumerate}
\end{thm}
Recall that the $\Aut^{loc}(M,\Tm)$-orbit of a point $x\in M$ is the set of points $y\in M$ for which there exists 
a local automorphism $f$ of $\Tm$, defined from a neighborhood of $x$ to a neighborhood of $y$,
and such that $f(x)=y$.

\subsection{Contact-Anosov flows}\label{C-A}
Theorem \ref{thm:GhyssansAnosov} was inspired and is a generalisation of Ghys' theorem \cite{G} classifying contact-Anosov flows 
with smooth invariant distributions
on compact three-manifolds 
(see section \ref{subsection-Anosovstrict} for definitions and more details about Anosov flows).  
Indeed, the one-parameter subgroup defined by a contact-Anosov flow is non-compact,
has dense orbits,
and preserves a strict path structure on a three-manifold.

\begin{thm}[\cite{G}]\label{thm:Ghys}
Let $(\varphi^t)$ be an Anosov flow on a compact connected three-dimensional manifold,
whose stable and unstable distributions $E^s$ and $E^u$ are $\cc^\infty$, 
and such that $E^s\oplus E^u$ is a contact distribution.
Then, up to finite coverings, $(\varphi^t)$ is $\cc^\infty$-orbitally equivalent to the geodesic flow of a compact hyperbolic surface.
\end{thm}

We recall that two flows are called \emph{orbitally equivalent}, if there exists a diffeomorphism between the manifolds supporting them,
sending the orbits of the first flow to the one of the second flow.

Ghys' Theorem was generalized in \cite{BFL} for Anosov flows in higher dimensions, using an adapted linear connection due to Kanai in \cite{K}, 
who obtained previous classification results in the case of geodesic flows.   Another important ingredient in this classification result is Gromov's theory of rigid geometric structures. 
On the other hand, a classification of partially hyperbolic diffeomorphisms preserving an \emph{enhanced path structure}
(that is a path structure with a fixed transverse direction to the contact distribution) on compact three manifolds was obtained in \cite{Mm}.  
{{An important distinction with the results of \cite{BFL,Mm} 
is that we don't impose in Theorem \ref{thm:GhyssansAnosov} the existence of special types of automorphisms, 
as in the case of Anosov flows or partially hyperbolic diffeomorphisms.}

\par On another note, the simplification of our method is apparent by the fact that we don't need to use Gromov's open dense orbit theorem (\cite{Gr}).
The main thrust is to use an adapted Cartan connection to obtain the two possible local models of the structure.  
Indeed, the existence of a  non-compact group of automorphisms  and a dense orbit easily imply that most components of the Cartan curvature vanish 
and that the structure is locally homogeneous.  

We think that a similar strategy could be fruitful to obtain new dynamical rigidity results,
{in higher dimensions, and for enhanced path structures.}

\subsection{Relations to other geometric structures}

Any strict path structure $\Tm=(E^1,E^2,\theta)$ on a three-dimensional manifold $M$ 
induces on $M$ a Lorentzian metric $g$ defined as follow: 
the Reeb vector field of $\theta$ is a unitary vector field for $g$, orthogonal to the plane $E^1\oplus E^2$, the directions $E^i$ are isotropic,
and $g(u,v)=d\theta(u,v)$ for any $(u,v)\in E^1\times E^2$. 
In other words, any strict path structure is a refinement of a three-dimensional Lorentzian metric.
As a consequence, our result follows from the more general classification
of compact Lorentzian three-manifolds having a non-compact isometry group, done by Frances in \cite{F1,F} (see also Zeghib \cite{Z}), 
which is a more difficult result.  Related to this theorem we should point out Fang's work \cite{Fa} classifying geometric (that is, preserving a  pseudo-Riemannian metric) Anosov flows in dimension five.
What is surprising, and seems interesting to us, is that
in the somehow specific case of strict path structures, 
the Cartan connection described in paragraph \ref{soussection:normalCartangeometrie} allows us to obtain the classification more easily.  

 In the case of CR-structures (see \cite{C,CM} for the construction of the appropriate Cartan connections), 
 if the group of automorphisms of a CR-manifold $M$ (with no compacity condition) 
 does not act properly, 
 then $M$ is the sphere or the Heisenberg group with the standard CR structures (\cite{Sc}). 
 Therefore, for a compact manifold, the group of automorphism is compact except for the sphere.  
 On the other hand for pseudo-hermitian structures (the analog of strict path structures in the CR context, see \cite{W})
 the automorphism group can be interpreted as the isometry group of a Riemannian metric,
 so that any compact pseudo-hermitian manifold has a compact automorphism group.
 \par Let us emphasize that in the case of path structures, the situation is more complicated.
 Indeed, 
 it is possible to construct compact flat path structures having a non-compact automorphism group, 
 but whose fundamental group is a free group. 
 Such examples have been constructed in  \cite{Mm2} as compactifications of 
 path structures defined on the unitary tangent bundle of non-compact hyperbolic surfaces.

\subsection*{Organization of the paper}
\par We organize the paper as follows. In section \ref{section:definition} we review the definitions of path structures, 
strict path structures, and define the models we work with.  
We also discuss the particular case of Anosov contact flows,
and explain in paragraph \ref{soussoussectionclassification} how Theorem \ref{thm:Ghys} follows
from Theorem \ref{thm:GhyssansAnosov}. 
The proof of Theorem \ref{thm:GhyssansAnosov} is  given in the next two sections.

In section \ref{section:fixedcontactform} we define a Cartan connection  associated to a strict path structure which is the main geometric tool used in the remaining section.

In the last section we prove Theorem \ref{thm:GhyssansAnosov}.  
We first prove in paragraph \ref{soussection:typecourbure} 
that certain components of the curvature must vanish if the automorphism group is large. 
This will imply that the manifold is locally homogeneous, 
and we identify in \ref{section:mutations} the two possible global models as one which is flat, and the other having constant curvature.  

\subsection*{Acknowledgements} 
We would like to thank Charles Frances for  several enlightening discussions.

\section{{Strict} path structures in dimension 3}\label{section:definition}

Path geometries are treated in detail in section 8.6 of \cite{IL} and in \cite{BGH} where the relation to second order differential equations is also explained.   

Le $M$ be a real three dimensional manifold and $\Fitan{M}$ be its tangent bundle.  

\begin{dfn}\begin{enumerate}
\item
A \emph{path structure} $\Lm=(E^1,E^2)$ on $M$ is a choice of two sub-bundles $E^1$ and $E^2$ in $\Fitan{M}$, such that 
$E^1 \cap E^2=\{ 0\}$ and $E^1 \oplus E^2$ is a contact distribution.  
\item A \emph{strict path structure} $\Tm=(E^1,E^2,\theta)$ on $M$ is a path structure with a fixed contact form $\theta$ such that $\ker\theta=E^1\oplus E^2$.  
\item A (local) \emph{automorphism} of $(M,\Tm)$ is a (local) diffeomorphism $f$ of $M$ that preserves $E^1$, $E^2$ and $\alpha$.
\end{enumerate}
\end{dfn}

The condition that $E^1 \oplus E^2$ be a contact distribution means that, locally,  there exists a one form $\theta$ on $M$ such that
$\ker \theta= E^1\oplus E^2$ and $\theta\wedge d\theta$ is never zero.  
On the other hand, for strict path structures we impose the existence of a globally defined contact form $\theta$.
Therefore, strict path structures are \emph{unimodular geometries}:
there exists a canonical volume form $\mu_{\Tm}=\theta\wedge d\theta$ on $M$, preserved by the automorphism group of $\Tm$ (in contrast, neither path structures nor enhanced path structures are unimodular).  There exists a unique vector field $R$ such that $d\theta(R,\cdot)=0$ and $\theta(R)=1$, called the \emph{Reeb vector field} of $\theta$,
that we will also call the Reeb vector field of the strict path structure $\Tm$.
In particular, the distribution $E^1\oplus E^2$ of a strict path structure $\Tm$ is thus oriented,
and the manifold $M$ supporting $\Tm$ is orientable.

Flat path structure   is the geometry of real flags in $\R^3$.  That is the geometry of the space of all couples $(p,l)$ where $p\in \R P^2$ and $l$ is a real projective line
containing $p$.  The space of flags is identified to the quotient
$$
\SL(3,\R)/B_\R
$$
where $B_\R$ is the Borel group of all real upper triangular matrices.  

\subsection{Flat strict path structure}\label{subsection-flatstrict}

\par The Heisenberg group 
\[
\mathbf{Heis}(3)
\coloneqq
\enstq{
\begin{pmatrix}
1 & 0 & 0 \\
y & 1 & 0 \\
z & x & 1 
\end{pmatrix}
}{(x,y,z)\in\R^3}
\]
is the model of strict path structures.  
We consider on $\mathbf{Heis}(3)$ the left-invariant structure 
$(\R\tilde{X},\R\tilde{Y},\widetilde{Z^*})$, where  $\tilde{X},\tilde{Y}$ are the left invariant vector fields and $\widetilde{Z^*}$ the left invariant 1-form 
induced by the basis
\begin{equation}\label{equation:baseheis3}
X=\left(\begin{smallmatrix}
0 & 0 & 0 \\
0 & 0 & 0 \\
0 & 1 & 0
\end{smallmatrix}\right),
Y=\left(\begin{smallmatrix}
0 & 0 & 0 \\
1 & 0 & 0 \\
0 & 0 & 0
\end{smallmatrix}\right),
Z=\left(\begin{smallmatrix}
0 & 0 & 0 \\
0 & 0 & 0 \\
1 & 0 & 0
\end{smallmatrix}\right)
\end{equation}
of its Lie algebra.
To describe the automorphism group of this structure, we introduce the 
subgroups
$$
                P=\enstq{
                \left ( \begin{array}{ccc}
                        a     &    0    &   0   \\
                       0    &    \frac{1}{a^2}    &   0\\
                       0     &  0 &     a

                \end{array} \right )\ \ \ 
                }{a\in\R^*}
                \subset 
G=\enstq{\left ( \begin{array}{ccc}

                        a     &    0    &   0   \\
                       y    &    \frac{1}{a^2}    &   0\\

                        z       &  x   &     a

                \end{array} \right )}{a\in\R^*,(x,y,z)\in\R^3}\ \ \                
$$
of $\SL(3,\R)$, and we decompose $G=\mathbf{Heis}(3)P$ through the diffeomorphism 
\[
\psi\colon(h,p)\in\mathbf{Heis}(3)\times P\mapsto hp\in G.
\]
More precisely, $\psi$ is a group isomorphism between $G$ and
the semi-direct product $\mathbf{Heis}(3)\rtimes P$, 
where $P$ acts on $\mathbf{Heis}(3)$ by $p\cdot h\coloneqq php^{-1}$.
We can thus identify $G$ with the subgroup $\mathbf{Heis}(3)\rtimes\mathcal{P}$
of affine group automorphisms of $\mathbf{Heis}(3)$ of the form $(h,\varphi)\colon x\in\mathbf{Heis}(3)\mapsto h\varphi(x)\in\mathbf{Heis}(3)$, 
where $h\in\mathbf{Heis}(3)$ and $\varphi\in\mathcal{P}$ is the conjugation by an element of $P$.
We use this identification to define the following left action of $G=\mathbf{Heis}(3)\rtimes P$ on $\mathbf{Heis}(3)$: 
\[
hp\cdot x\coloneqq h(pxp^{-1})
\]
for any $(h,p)\in\mathbf{Heis}(3)\times P$ and $x\in\mathbf{Heis}(3)$.
\par This action being transitive, it induces an identification of $\mathbf{Heis}(3)$ with $G/P$,
by choosing the identity $e$ for base-point.
It is easy to verify that $\R X$, $\R Y$ and $Z$ are fixed by the adjoint action of $P$,
so that $G$ acts on $\mathbf{Heis}(3)$ by automorphisms of its strict path structure
($G$ is actually the whole automorphism group of this structure).
If $\Gamma$ is a discrete subgroup of $G$ acting freely, properly and cocompactly on $\mathbf{Heis}(3)$, 
we will always implicitly endow the quotient $\Gamma\backslash\mathbf{Heis}(3)$ with the induced strict path structure.
\par It is possible to construct such compact quotients 
admitting a non-compact automorphism group (for the compact-open topology).
Indeed, with $\mathbf{Heis}_\Z(3)$ the cocompact lattice of elements of $\mathbf{Heis}(3)$ having integer entries, 
it is easy to find a group automorphism $\phi_0$ of $\mathbf{Heis}(3)$ preserving $\mathbf{Heis}_\Z(3)$, 
and such that $\Diff{e}{\phi_0}$ is diagonalizable with eigenvalues $(\lambda,\lambda^{-1},1)$, 
where $\abs{\lambda}<1$ and $\R Z$ is the eigenline for the eigenvalue $1$.
Now, there exists an automorphism $\Phi$ of $\mathbf{Heis}(3)$ such that 
$\phi\coloneqq\Phi\circ\phi_0\circ\Phi^{-1}\in\mathcal{P}$, 
and since $\phi$ preserves the cocompact lattice $\Gamma\coloneqq\Phi(\mathbf{Heis}_\Z(3))$,
it induces a diffeomorphism $\bar{\phi}$ of the compact quotient $\Gamma\backslash\mathbf{Heis}(3)$.
The dynamic of $\bar{\phi}$ is similar to the time-one of an Anosov flow, 
in the sense that it contracts one of the two directions $(\R\tilde{X},\R\tilde{Y})$ and expands the other one.
In particular, $\bar{\phi}$ generates a non relatively compact subgroup of automorphisms of $\Gamma\backslash\mathbf{Heis}(3)$. 
\par The diffeomorphisms having this type of dynamics are called \emph{partially hyperbolic},
and the purpose of \cite{Mm} is precisely to classify those partially hyperbolic diffeomorphisms 
of contact type whose invariant distributions are smooth, which is the discrete-time counterpart of Ghys Theorem \ref{thm:Ghys}.

\subsection{$\widetilde{\mathbf{SL}}(2,\R)$ as a strict path structure}\label{subsection-SL2strict}
\label{subsubsection:deformationflotsgeodesiques}

Consider the universal cover $\pi\colon\widetilde{\SL}(2,\R)\to\PSL(2,\R)$ of $\PSL(2,\R)$.
Let us use the following usual basis for its Lie algebra $\mathfrak{sl}_2$:
\begin{equation}\label{equation:basecanoniquesl2}
E=
\left(\begin{smallmatrix}
         0 & 1 \\
         0 & 0 
        \end{smallmatrix}\right),
         F=\left(\begin{smallmatrix}
         0 & 0 \\
         1 & 0 
        \end{smallmatrix}\right), \text{~and~}
         H=\left(\begin{smallmatrix}
         1 & 0 \\
         0 & -1 
        \end{smallmatrix}\right).
\end{equation}
The Lie bracket relation $[E,F]=H$ shows that this basis defines
a left-invariant strict path structure
$(\R\tilde{E},\R\tilde{F},\widetilde{H^*})$ on $\widetilde{\SL}(2,\R)$, 
where $H^*\in(\mathfrak{sl}_2)^*$ denotes the linear form dual to $H$ ($H^*(H)=1,H^*(E)=H^*(F)=0$).
In order to describe the automorphism group of this strict path structure,
we define the \emph{right diagonal flow} of $\PSL(2,\R)$
as the flow by right translations of the one-parameter subgroup $A\coloneqq\{a^t\}_{t\in\R}\subset \PSL(2,\R)$, where
\begin{equation}\label{equation:at}
a^t\coloneqq
\begin{bmatrix}
e^t & 0 \\
0 & e^{-t}
\end{bmatrix}.
\end{equation}
We denote by 
$\tilde{A}$ the 
subgroup of $\widetilde{\SL}(2,\R)$ which projects to $A$ on $\PSL(2,\R)$, 
by $\{\tilde{a}^t\}$ the one-parameter subgroup inside $\tilde{A}$ which projects to $\{a^t\}$,
and we continue to call \emph{right diagonal flow} its flow by right translations on $\widetilde{\SL}(2,\R)$.

We define a (left) action of $\widetilde{\SL}(2,\R)\times\tilde{A}$ on $\widetilde{\SL}(2,\R)$ by $(g,a)\cdot x\coloneqq gxa$.
Since the action of the right diagonal flow preserves the strict path structure
(because the adjoint action of $H$ preserves $\R E$ and $\R F$),
$\widetilde{\SL}(2,\R)\times\tilde{A}$ is contained in the automorphism group of $\widetilde{\SL}(2,\R)$
(it is actually not difficult to verify that this is its whole automorphism group).

Let $\Gamma_0$ be a cocompact lattice of $\widetilde{\SL}(2,\R)$,
$\rho\colon\Gamma_0\to\tilde{A}$ be a morphism, and let us denote by $gr(\rho)$ the graph of $\rho$,
which is a discrete subgroup of $\widetilde{\SL}(2,\R)\times\tilde{A}$.
If $gr(\rho)$ acts freely, properly and cocompactly on $\widetilde{\SL}(2,\R)$ then following \cite{S}, 
we will say that $\rho$ is \emph{admissible}.
The subgroups $\Gamma$ of $\widetilde{\SL}(2,\R)\times\tilde{A}$ of the form $gr(\rho)$ with $\rho$ an admissible morphism from a cocompact lattice,
will be called \emph{admissible discrete graph-groups}.
Consequently, any quotient $\Gamma\backslash\widetilde{\SL}(2,\R)$ by an admissible discrete graph-group of $\widetilde{\SL}(2,\R)\times\tilde{A}$
inherits of the structure of $\widetilde{\SL}(2,\R)$, and we will always endow such a quotient $\Gamma\backslash\widetilde{\SL}(2,\R)$
with this induced strict path structure.


\subsection{Contact Anosov Flows}\label{subsection-Anosovstrict}
It happens that
strict path structures $\Tm=(E^1,E^2,\theta)$ naturally appear in an important dynamical situation: 
the one of a compact three-manifold $M$ endowed with a \emph{contact Anosov flow} with smooth invariant distributions.
\subsubsection{Anosov flows}\label{soussoussectionAnosov}
We first recall the definition of an Anosov flow.
\begin{dfn}\label{dfn:Anosovflow}
A non-singular flow $(\varphi^t)$ of class $\cc^\infty$ of a compact manifold $M$ is called \emph{Anosov},
if its differential preserves a splitting $\Fitan{M}=E^s\oplus E^0\oplus E^u$ of the tangent bundle, where $E^0=\R X^0$
with $X^0$ the (non-singular) vector field generating $(\varphi^t)$,
and where $E^s$ and $E^u$ are non-trivial distributions
verifying the following estimates (with respect to any Riemannian metric on $M$).
\begin{enumerate}
\item The \emph{stable distribution} $E^s$ is \emph{uniformly contracted} by $(\varphi^t)$, 
\emph{i.e.} there are two constants $C>0$ and $0<\lambda<1$ such that
for any $t\in\R$ and $x\in M$: 
\begin{equation}\label{equationcontractionuniforme}
\norme{\Diff{x}{\varphi^t}\restreinta_{E^s}}\leq C\lambda^t.
\end{equation}
\item The \emph{unstable distribution} $E^s$ is \emph{uniformly expanded} by $(\varphi^t)$, \emph{i.e.} uniformly contracted by $(\varphi^{-t})$.
\end{enumerate}
\end{dfn}

\par The geodesic flow of a compact hyperbolic surface $\Sigma$, acting on its unitary tangent bundle $\Fiunitan{\Sigma}$,
is an important example of three-dimensional Anosov flow 
(it is in fact the historical motivation for the study of these flows, see \cite{A}).
To have a geometric image in mind, the projections on $\Sigma$ of the stable and unstable foliations of this flow are projections
of horocycles of the hyperbolic plane $\Hn{2}$.
These examples also have a natural description in terms of homogeneous spaces.
Let $\Sigma$ be the quotient of $\Hn{2}$ by a cocompact lattice $\Gamma_0$ (without torsion)
of its group of orientation-preserving isometries, identified with $\PSL(2,\R)$.
Since $\PSL(2,\R)$ acts simply transitively on the unitary tangent bundle of $\Hn{2}$,
one verifies that the geodesic flow on $\Fiunitan{\Sigma}$ is smoothly conjugated to 
(a constant time-change of) the right diagonal flow on the quotient $\Gamma_0\backslash\PSL(2,\R)$,
that was defined in the paragraph \ref{subsection-SL2strict}.

\par Let us emphasize that in the definition of Anosov flows, \emph{no regularity} 
is requested on the stable and unstable distributions.
Even if they are automatically H\"older continuous (thanks to the estimates \eqref{equationcontractionuniforme}),
$E^s$ and $E^u$ have in general no reason to be differentiable (even if the flow is $\cc^\infty$).
This ``lack of regularity'' is in fact much more than a detail, as it is in a sense responsible for the abundance of Anosov flows. 
For instance, the use of surgery methods allowed Handel and Thurston to prove in \cite{HT} the existence of \emph{non-algebraic Anosov flows},
which are not conjugated to any Anosov flow with smooth distributions.
\par From that perspective, geodesic flows of compact hyperbolic surfaces are very specific: 
their stable and unstable distributions both are $\cc^\infty$
(because they arise from left-invariant distributions on the Lie group $\PSL(2,\R)$). 
Their sum $E^s\oplus E^u$ is moreover a contact distribution, 
and Anosov flows verifying this last property are called \emph{contact-Anosov}.
If $(\varphi^t)$ is a contact-Anosov flow with smooth invariant distributions,
and $X^0$ its infinitesimal generator, then we define the \emph{canonical one-form} $\theta$ of $(\varphi^t)$ by 
$\theta\restreinta_{E^s\oplus E^u}=0$ and $\theta(X^0)=1$.
This is a contact form of kernel $E^s\oplus E^u$,
and by construction, $(\varphi^t)$ preserves the strict path structure $\Tm=(E^s,E^u,\theta)$
that we will call its \emph{canonical} structure.
Note that the structures obtained in this way have a purely geometrical specificity: 
the Reeb flow of their contact form is a flow of automorphisms of the structure $\Tm$
(which has no reason to be true in general).
Indeed, the Reeb vector field of the canonical structure of $(\varphi^t)$ is its generator $X^0$, 
so that the Anosov flow itself is encoded in the structure $\Tm$. 

\subsubsection{Deformations of geodesic flows}\label{soussoussectiondeformations}
We now use the structures described in the paragraph \ref{subsection-SL2strict}
to obtain finite coverings and time changes of geodesic flows of compact hyperbolic surfaces.
\par Let $\Gamma_0$ be a cocompact lattice of $\widetilde{\SL}(2,\R)$.
The \emph{finiteness of the level} proved by Salein in \cite[Th\'eor\`eme 3.3.2.3]{S}
(see also \cite{KR}) shows that its projection $\bar{\Gamma}_0=\pi(\Gamma_0)$ in $\PSL(2,\R)$ is 
a cocompact lattice of $\PSL(2,\R)$, and that 
the projection $\Gamma_0\backslash\widetilde{\SL}(2,\R)\to\bar{\Gamma}_0\backslash\PSL(2,\R)$ induced by $\pi$ is a finite covering.
According to Selberg's Lemma, $\bar{\Gamma}_0$ moreover contains a finite index subgroup $\Gamma_0'\subset\PSL(2,\R)$ 
without torsion, and the right diagonal flow on $\Gamma_0\backslash\widetilde{\SL}(2,\R)$ is finally a finite covering 
of the geodesic flow of the hyperbolic surface $\Gamma_0'\backslash\Hn{2}$.
As such, it is in particular an Anosov flow.
\par Now if 
$\rho\colon\Gamma_0\to\tilde{A}$ is an admissible morphism and 
$\Gamma=gr(\rho)\subset\widetilde{\SL}(2,\R)\times\tilde{A}$ 
its graph-group,
the right diagonal flow of $\widetilde{\SL}(2,\R)$ induces a flow on the compact quotient $\Gamma\backslash\widetilde{\SL}(2,\R)$
(because $\tilde{A}$ is abelian),
and this flow remains Anosov.\footnote{This is non-trivial, at least for us, and can be explained as follows:
\cite[Prop. 4.2]{Z} shows that this flow is quasi-Anosov, 
and in dimension three, quasi-Anosov flows are Anosov according to \cite[Theorem A]{M}.}
In fact, it is possible to prove that $\Gamma\backslash\widetilde{\SL}(2,\R)$ 
is diffeomorphic to $\Gamma_0\backslash\widetilde{\SL}(2,\R)$,
and that we can choose a diffeomorphism sending orbits of the right diagonal flow 
to orbits of the right diagonal flow 
(this is a consequence of a result of Haefliger in \cite{H}, and the reader can find more details about this in \cite[Th\'eor\`eme 6.5]{G}).
\par We have described a family of
contact-Anosov flows that are orbit-equivalent to finite covers of geodesic flows,
and whose invariant distributions are smooth (since they arise from left-invariant distributions of $\widetilde{\SL}(2,\R)$). 

\subsubsection{Classification in dimension three}\label{soussoussectionclassification}
\par On the other hand, flows on the Heisenberg group are more restricted. One can ask if  
 it is possible for the Reeb flow of a compact quotient $\Gamma\backslash\mathbf{Heis}(3)$ to be \emph{Anosov} ?
The following result answers by the negative.
This statement is probably already known, but we did not find it in the literature, 
so we suggest a proof in the appendix \ref{subsection:pasHeisenberg}.
\begin{prop}\label{po:pasHeis}
Let $\Gamma$ be a discrete subgroup of $G$ acting freely, properly and cocompactly on $\mathbf{Heis}(3)$.
Then on the quotient $\Gamma\backslash\mathbf{Heis}(3)$, 
the central flow of $\mathbf{Heis}(3)$
acts periodically.
\end{prop}
Let us emphasize that the flow by left translations of the center of $\mathbf{Heis}(3)$,
that we call its \emph{central flow}, induces a flow on any quotient $\Gamma\backslash\mathbf{Heis}(3)$,
which is precisely its Reeb flow.  
\par With the help of Proposition \ref{po:pasHeis}, 
Theorem \ref{thm:GhyssansAnosov} implies that the contact Anosov flows with smooth invariant distributions
previously described in the paragraph \ref{subsubsection:deformationflotsgeodesiques} are in fact the only possible examples,
which implies in particular Ghys Theorem \ref{thm:Ghys}.
\begin{cor}[\cite{G}]\label{thm:Ghysprecision}
Let $(\varphi^t)$ be a contact-Anosov flow with smooth invariant distributions 
on a compact connected three-dimensional manifold.
Then $(\varphi^t)$ is smoothly conjugated to a constant time-change of
the right diagonal flow on a compact quotient $\Gamma\backslash\widetilde{\SL}(2,\R)$,
with $\Gamma$ an admissible discrete graph-group of $\widetilde{\SL}(2,\R)\times\tilde{A}$.
\end{cor}
\Pf
We apply Theorem \ref{thm:GhyssansAnosov} to the canonical strict path structure $\Tm$ of 
a contact-Anosov flow $(\varphi^t)$ with smooth invariant distributions.
Since $(\varphi^t)$ preserves a contact form and hence a volume, it has no wandering points and is thus topologically transitive
(see \cite{Pu}). 
In particular, the automorphism group of $\Tm$ has a dense orbit.
therefore, $\Tm$ verifies the hypotheses of Theorem \ref{thm:GhyssansAnosov}, 
and since $(\varphi^t)$ is the Reeb flow of $\Tm$ and is not periodic, Proposition \ref{po:pasHeis} says that
$\Tm$ is isomorphic to $\Gamma\backslash\widetilde{\SL}(2,\R)$ 
with $\Gamma$ a discrete subgroup of $\widetilde{\SL}(2,\R)\times\tilde{A}$
acting freely, properly and cocompactly on $\widetilde{\SL}(2,\R)$.
According to the work of Kulkarni-Raymond in \cite{KR} and Salein in \cite[Theorem 3.3.2.3]{S} (see also Tholozan, \cite[Lemma 4.3.1]{Tho}),
such a group $\Gamma$ is the graph of an admissible morphism
$\rho$ from a cocompact lattice $\Gamma_0$ of $\widetilde{\SL}(2,\R)$
(details can be found in \cite[Fact 8.1]{Mm}).
This concludes the proof of Corollary \ref{thm:Ghysprecision}.
\EPf

\section{A Cartan connection for strict path structures} \label{section:fixedcontactform}

Cartan connections were used by Cartan in his numerous studies and classifications of geometric structures.
The main idea is to approximate a geometric structure by the model space $G/P$, 
 using a $P$-principal bundle with a Cartan connection (see the definition below) which mimics the Maurer-Cartan form of $G$ in the case of the model.  The construction of the appropriate Cartan connection associated to a geometric structure can sometimes be difficult.

A Cartan connection for path structures was first obtained by Cartan in \cite{Ca}.  One can refer to \cite{IL, BGH} for a modern treatment.   
 In this section we will describe a Cartan connection for  strict path structures.  In this reduction of the structure, the
 construction of a Cartan connection is much easier. 
 It was previously used in \cite{FV} in a slightly more general context of complexified tangent spaces of real three dimensional manifolds.  
 In that case we named a path structure with a fixed contact form a pseudo flag structure, in analogy to pseudo-hermitian structures.  
 
 \subsection{Preliminaries on Cartan geometries}
 We first define the main notions that we will use about Cartan geometries,
 and recall general facts.
 The reader will find 
 more details, and proofs of all the claims given concerning Cartan geometries,
 in \cite{Sh} for instance.
 \par One starts with a model space $X=G/P$ where $P$ is a closed subgroup of a Lie group $G$.  
 We see then $G$ as a right $P$-principal bundle over $X$, and the Maurer-Cartan form
 $\varpi : \Fitan{G}\rightarrow \mathfrak{g}$ can be viewed as an identification of the tangent bundle of $G$ with its Lie algebra $\mathfrak{g}$, 
 through left invariant vector fields.
Denoting by $R_p$ the right action of an element $p\in P$,
$\varpi$ satisfies the equivariance property 
 $ R_h^*\varpi= Ad_{h^{-1}}\varpi$.
 
 The approximation by the model space of a geometric structure is described using a Cartan connection on an appropriate principal bundle:
 
 \begin{dfn} A \emph{Cartan geometry} $\mathcal{C}=(\hat{M},\varpi)$ modelled on $X$ over a manifold $M$ is a right $P$-principal bundle equipped with a form 
 $\varpi : \Fitan{\hat{M}}\rightarrow \mathfrak{g}$ (called a Cartan connection) satisfying the following properties.
 \begin{enumerate}
 \item At each point $\hat{x}\in \hat{M}$, $\varpi_{\hat{x}} : \Tan{\hat{x}}{\hat{M}}\rightarrow \mathfrak{g}$  is an isomorphism.
 \item If $X\in {\mathfrak p}$, and if $X^*\colon\hat{x}\in\hat{M}\mapsto \frac{d}{dt}\restreinta_{t=0}\hat{x}\cdot \exp(t X)$ 
 denotes the fundamental vector field associated to $X$, then
	$\varpi(X^*)=X$.
 \item If $p\in P$ then  $R_p^* \varpi=Ad_{p^{-1}}\varpi$.
 \end{enumerate}
 \end{dfn}

 A \emph{(local) automorphism} $f: M\rightarrow M$ of the Cartan geometry $\mathcal{C}$ 
 is a (local) diffeomorphism $f$ of $M$ that lifts to 
 a $P$-equivariant (local) diffeomorphism $\hat{f } : \hat{M}\rightarrow \hat{M}$ such that $\hat{f } ^*\varpi=\varpi$.
 Such a lift to $\hat{M}$ of an automorphism $f$ is unique, and will always be denoted by $\hat{f}$.
 
 The \emph{curvature form} of a Cartan geometry $(\hat{M},\varpi)$ over $M$ is the $\mathfrak{g}$-valued two-form defined on $\hat{M}$ by
 $$
 \Omega(X,Y)= d\varpi(X,Y)+[\varpi(X),\varpi(Y)].
 $$
Observe that in the flat case $G\rightarrow G/P$, the Maurer-Cartan equation implies that $\Omega=0$.  Also, due to its construction, the curvature form vanishes on vertical vectors.  One may use the connection form to define the curvature map
$K:\hat{M}\rightarrow L(\Lambda^2\mathfrak{g},\mathfrak{g})$ by
\begin{equation}\label{equation:courbureT}
K(\hat{x})(u,v)=\Omega_{\hat{x}}(\varpi_{\hat{x}}^{-1}(u),\varpi_{\hat{x}}^{-1}(u)).
\end{equation}
The vanishing of $\Omega$ on vertical vectors implies that the curvature map induces a map 
$K:\hat{M}\rightarrow L(\Lambda^2(\mathfrak{g}/\mathfrak{p}),\mathfrak{g})$  which we denote by the same symbol.   
Observe also that an automorphism of $M$ preserves $K$, that is $K\circ \hat{f}=K$.  

The curvature form satisfies the following equivariance property
$$
 R^*_p\Omega =Ad_{p^{-1}}\Omega.
 $$
 This equivariance  translates into an equivariance of the curvature map under the action of $P$,
 where the right action of $P$ on $L(\Lambda^2(\mathfrak{g}/\mathfrak{p}),\mathfrak{g})$ is given by
 $$
 (K.p)(u,v)=Ad_{p^{-1}}(K(Ad_{p}u,Ad_{p}v)).
 $$
 
\subsection{Normal Cartan geometry of a strict path structure}\label{soussection:normalCartangeometrie}
We now go back to strict path structures, 
by considering the specific case of Cartan geometries modelled on $\mathbf{Heis}(3)$, the flat model of strict path structures
 introduced in paragraph \ref{subsection-flatstrict}.
 So $G$ denotes from now on the subgroup of $\SL(3,\R)$ defined by
 $$
 G=\enstq{\left ( \begin{array}{ccc}

                        a     &    0    &   0   \\
                       y    &    \frac{1}{a^2}    &   0\\

                        z       &  x &    a

                \end{array} \right )}{a\in\R^*,(x,y,z)\in\R^3}
$$ 
and $P\subset G$ the subgroup of $G$ defined by
$$
 P=\left\{\left ( \begin{array}{ccc}

                        a     &    0    &   0   \\
                       0    &    \frac{1}{a^2}    &   0\\

                        0      &     0 &     a

                \end{array} \right )\right\}.\ \ \
$$
 
 Writing the Maurer-Cartan form of $G$ as the matrix
$$
\left ( \begin{array}{ccc}

                        w     &    0    &   0   \\
                       \omega^2    &    -2w   &   0\\

                       \omega      &  \omega^1   &     w

                \end{array} \right )
$$
one obtains the Maurer-Cartan equations:
$$d \omega-\omega^2\wedge \omega^1=0$$
$$d\omega^1 +3w\wedge \omega^1=0$$
$$d\omega^2 -3w\wedge \omega^2=0$$
$$d w=0.$$

{
We saw in the paragraph \ref{subsection-flatstrict} that $G$ is the automorphism group of the canonical left-invariant strict path structure of $\mathbf{Heis}(3)$,
and that its action induces an identification of $\mathbf{Heis}(3)$ with the homogeneous space $X=G/P$.
}

Let $M$ be a three-manifold equipped with a strict path structure $\Tm=(E^1,E^2,\theta)$ of Reeb vector field $R$.
Now let  $X_1\in E^1$, $X_2\in E^2$ be such that $d\theta(X_1,X_2)=1$.  
The dual coframe of $(X^1,X^2,R)$ is $(\theta^1,\theta^2,\theta)$,
for two 1-forms $\theta_1$ and $\theta_2$ verifying $d\theta=\theta^1\wedge \theta^2$.   

At any point $x\in M$,
one can look at the coframes of the form
$$
\omega^1=\frac{1}{a}\theta^1(x), \  \omega^2=a\theta^2(x), \ \omega=\theta(x)
$$
for $a\in\R^*$.
\begin{dfn}
We denote by $\pi: \hat{M}\rightarrow M$ the  {$P$}-coframe bundle over $M$ 
  given by the set of coframes  $(\omega, \omega^{1}, \omega^{2})$ of the above form.
The structure group $P$ acts on $\hat{M}$ as follows
$$
(\omega', \omega'^{1},
  \omega'^{2})=
(\omega, \omega^{1},
  \omega^{2})
 \left ( \begin{array}{cccc}

                        1       &    0               	&       0 \\

                        0  	&      \frac{1}{a} 		& 0             \\
                        0   	&  0 			&  a  \\

                \end{array} \right )
$$
where   $a \in \R^*$.
\end{dfn}

We will denote the tautological forms defined by $\omega^1,\omega^2,\omega$ using the same letters. 
That is, we write $\omega^i$ at the coframe $(\omega^1,\omega^2,\omega)$ to be $\pi^*(\omega^i)$.

\begin{prop}\label{theta11}
There exist unique forms $w$, $\tau^1$ and $\tau^2$ on $\hat{M}$ such that 
\begin{equation}\label{dtheta1}
d\omega^1=3\omega^1\wedge w+\omega\wedge\tau^1,
\end{equation}
\begin{equation}\label{dtheta2}
d\omega^2=-3\omega^2\wedge w+\omega\wedge\tau^2,
\end{equation}
\begin{equation}
{
\text{and~}\tau^1\wedge \omega^2=\tau^2\wedge \omega^1=0.
}
\end{equation}
{Moreover, $w=\frac{da}{3a}  \mod \omega^1,\ \omega^{2},\ \omega$.}
\end{prop}

Observe that the condition $\tau^1\wedge \omega^2=\tau^2\wedge \omega^1=0$ implies that we may write
$\tau^1=\tau^1_2\omega^2$ and $\tau^2=\tau^2_1\omega^1$.

\Pf  
 First observe that, writing $d\theta^1= A^1_{10}\theta^1\wedge \theta +A^1_{20}\theta^2\wedge \theta+A^1_{12}\theta^1\wedge \theta^2$ and 
$d\theta^2= A^2_{10}\theta^1\wedge \theta +A^2_{20}\theta^2\wedge \theta+A^2_{12}\theta^1\wedge \theta^2$,
we obtain, differentiating  $d\theta=-\theta^1\wedge \theta^2$, that is $0=d\theta^1\wedge \theta^2-\theta^1\wedge d\theta^2$, that $A^1_{10}=-A^2_{20}$.

We differentiate the tautological forms $\omega^1=\frac{1}{a}\theta^1$ and $\omega^2={a}\theta^2$ (by abuse of notation we write $\theta^i=\pi^*(\theta^i)$).  We have
$$
d\omega^1= d(\frac{1}{a})\wedge \theta^1+\frac{1}{a}d\theta^1=-\frac{da}{a}\wedge\omega^1+\frac{1}{a}(A^1_{10}\theta^1\wedge \theta +A^1_{20}\theta^2\wedge \theta+A^1_{12}\theta^1\wedge \theta^2)
$$
$$
=(-\frac{da}{a}-A^1_{10}\omega- aA^1_{12}\omega^2 )\wedge\omega^1+\frac{1}{a^2}A^1_{20}\omega^2\wedge \omega
$$
which can be written
$$
d\omega^1= -3w\wedge\omega^1+\omega\wedge \tau^1,
$$
for a horizontal form $\tau^1$.  We can choose $\tau^1$ therefore satisfying $\tau^1\wedge \omega^2=0$.
Analogously,
$$
d\omega^2= da\wedge \theta^2+ad\theta^2=\frac{da}{a}\wedge\omega^2+{a}(A^2_{10}\theta^1\wedge \theta +A^2_{20}\theta^2\wedge \theta+A^2_{12}\theta^1\wedge \theta^2)
$$
$$
=(\frac{da}{a}-A^2_{20}\omega+\frac{1}{a}A^2_{12}\omega^1 )\wedge\omega^2+{a^2}A^2_{10}\omega^1\wedge \omega
$$
We write therefore 
$$
d\omega^2= 3w \wedge\omega^2+\omega\wedge \tau^2,
$$
for a horizontal form $\tau^2$ satisfying  $\tau^2\wedge \omega^1=0$.  Remark that 
$w= \frac{da}{a}-A^2_{20}\omega+\frac{1}{a}A^2_{12}\omega^1+ aA^1_{12}\omega^2$ (because $A^1_{10}=-A^2_{20}$).

If other forms $w', \tau'^1, \tau'^2$ with the condition  $\tau'^1\wedge \omega^2=\tau'^2\wedge \omega^1=0$ satisfy these equations we obtain
$$
3(w-w')\wedge \omega^1+ \omega\wedge (\tau^1-\tau'^1)=0
$$
and
$$
3(w-w')\wedge \omega^2+ \omega\wedge (\tau^2-\tau'^2)=0
$$
It is clear then that these equations imply that $w=w', \tau^1=\tau'^1$ and $\tau^2= \tau'^2$.
\EPf

The Cartan connection can then be written as
$$
\varpi=\left ( \begin{array}{ccc}

                        w     &    0    &   0   \\
                       \omega^2    &    -2w   &   0\\

                       \omega      &  \omega^1   &     w

                \end{array} \right )
$$
satisfying  for an element 
$$
p=
\left( \begin{array}{ccc}

                        a  	&    	0     	&       0 \\

                       0				 &     \frac{1}{a^2}		&		0\\
                       0   				&  	0 				&		a
                      
       \end{array} \right )\in P,
$$
$$
R^*_p\varpi=Ad_{p^{-1}}\varpi=
\left( \begin{array}{ccc}

                        w  	&    	0     	&       0 \\

                       a^{3}\omega^2			 &     -2w		&		0\\
                       \omega   				&  	a^{-3}\omega^1				&		w
                      
       \end{array} \right ).
$$
The two other properties defining a Cartan connection being easily verified,
$\mathcal{C}=(\hat{M},\varpi)$ is finally a Cartan geometry modelled on $G/P$ over $M$, 
which will be called the \emph{normal Cartan geometry} of $(M,\Tm)$.
{
Note that the strict path structure $\Tm=(E^1,E^2,\theta)$ of Reeb vector field $R$ is determined by its normal Cartan geometry $\mathcal{C}$
through the following relations:
\begin{equation}\label{equationrelationCT}
E^1=\Diff{}{\pi}\circ\varpi^{-1}(\R X), E^2=\Diff{}{\pi}\circ\varpi^{-1}(\R Y), R=\Diff{}{\pi}\circ\varpi^{-1}(Z).   
\end{equation}
}
The fundamental property of the normal Cartan geometry, is that it has the same local automorphisms as $\Tm$.
\begin{thm}\label{theorem:geomCartannormaleT}
Let $(M,\Tm)$ be a three-dimensional strict path structure, and $\mathcal{C}$ be the normal Cartan geometry of $\Tm$.
For $f$ a local diffeomorphism of $M$, $f$ is an automorphism of $\Tm$ if, and only if it is an automorphism of $\mathcal{C}$.
\end{thm}
Note that any automorphism $f$ of $\Tm$ acts on $\hat{M}$ by the following relation: for 
$\hat{x}=(\omega, \omega^{1}, \omega^{2})\in\pi^{-1}(x)$, 
$\hat{f}(\hat{x})$ is the coframe $(\omega, \omega^{1}, \omega^{2})\circ\Diff{f(x)}{f^{-1}}$ at $f(x)$.
In other words, if $s\colon U\to \hat{M}$ is a local section,
$s(x)=(\theta,\theta^1,\theta^2)_x$ with $(\theta,\theta^1,\theta^2)$
a local coframe field,
then:
\begin{equation}\label{equation:actionautomorphismes}
 \hat{f}(s(x))=({f^{-1}}^*s)_{f(x)}.
\end{equation}

The curvature form of the normal Cartan geometry of $\Tm$ is given by
\begin{equation*}\label{equationcurvatureform}
\Omega=\left ( \begin{array}{ccc}
                      dw&   0      &  0   \\

                        \tau^2_1 \omega\wedge \omega^1 &   -2dw    &  0  \\

                           0  & \tau^1_2 \omega\wedge \omega^2 &    dw
               \end{array} \right ).
\end{equation*}
Using the notation \eqref{equation:baseheis3} for the basis $(X,Y,Z)$ of $\mathfrak{heis}(3)$,
this means that the curvature map $K$ of $\mathcal{C}$ has values in the following subspace of  
$L(\Lambda^2(\mathfrak{g}/\mathfrak{p}),\mathfrak{g})$:
\begin{equation}\label{equationW}
\mathcal{W}=
\enstq{K}{
K(\bar{Z},\bar{X})=
\left(\begin{smallmatrix} 
                      W^1 & 0 & 0 \\
 					\tau^2_1 & -2W^1 & 0 \\
 					0 & 0 & W^1
 			        \end{smallmatrix}\right),
K(\bar{Z},\bar{Y})=         \left(\begin{smallmatrix} 
                      W^2 & 0 & 0 \\
 				    0 & -2W^2 & 0 \\
 				    0 & \tau^1_2 & W^2
 			        \end{smallmatrix}\right),                    
K(\bar{X},\bar{Y})=
                      \left(\begin{smallmatrix} 
                      R & 0 & 0 \\
 					0 & -2R & 0 \\
 					0 & 0 & R
 			        \end{smallmatrix}\right)
}.
\end{equation}

\begin{lemma}\label{lemma:ReebKilling}
 If $W^1$, $W^2$, $\tau^2_1$ and $\tau^1_2$ identically vanish,
 then the Reeb vector field is a Killing field of the strict path structure.
\end{lemma}
\Pf
Indeed in this case,
denoting $\tilde{Z}\coloneqq\varpi^{-1}(Z)$ we have for any
$u\in\Fitan{\hat{M}}$:
$d\varpi(\tilde{Z},u)=K(Z,\varpi(u))-[Z,\varpi(u)]=0$
since $Z$ is in the center of $\mathfrak{g}$, and thus:
 $\mathcal{L}_{\tilde{Z}}\varpi=
 d(\varpi(\tilde{Z}))+d\varpi(\tilde{Z},\cdot)=0$,
 which shows that $\tilde{Z}$ is the lift of a Killing vector field
 of $\mathcal{C}$.
 The projection of $\tilde{Z}$ being the Reeb vector field,
 this proves the claim.
\EPf

A direct calculation shows that $\mathcal{W}$ is $P$-invariant, and if the components of $K\in\mathcal{W}$ 
are denoted by $W^1,\tau^2_1,W^2,\tau^1_2,R$ with the 
above notations, 
and 
\[
p=
\begin{pmatrix}
a & 0 & 0 \\
0 & a^{-2} & 0 \\
0 & 0 & a
\end{pmatrix},
\]
then $K\cdot p$ has coordinates $(\tilde{W}^1,\tilde{\tau}^2_1,\tilde{W}^2,\tilde{\tau}^1_2,\tilde{R})$ with
\begin{equation}\label{equationactionAWT}
\tilde{W}^1=a^3W^1,\tilde{\tau}^2_1=a^6\tau^2_1,\tilde{W}^2=a^{-3}W^2,\tilde{\tau}^1_2=a^{-6}\tau^1_2,\tilde{R}=R.
\end{equation}

\vspace{.5cm}

{\bf Example:}  Let $\SU(2)$ with its left invariant vector fields defined by a Lie algebra basis $(X^1,X^2,Y)$ with 
$[X^1,X^2]=Y$ and cyclic permutation of this commutation relation.   The Maurer-Cartan equations, for a dual basis  $\theta^1,\theta^2,\theta$, are:

$$d \theta+\theta^1\wedge \theta^2=0,$$
$$d\theta^1 +\theta^2 \wedge \theta=0,$$
$$
d\theta^2 +\theta \wedge \theta^1=0.
$$

Consider now  the strict path structure defined by $(\R X^1,\R X^2,Y)$.
Defining the tautological forms $\omega^1=\frac{1}{a}\theta^1$ and $\omega^2={a}\theta^2$, computng derivatives and comparing with Proposition \ref{theta11} we obtain that, 
$$
d\omega^1=3\omega^1\wedge w+\omega\wedge\tau^1,
$$
with $3w=\frac{da}{a}$ and $\tau^1= -\frac{1}{a^2}\omega^2$.
Analogously
$$
d\omega^2=-3\omega^2\wedge w+\omega\wedge\tau^2,
$$
with $\tau^2= {a^2}\omega^1$.

The curvature form of the connection is given by
$$
\Omega=\left ( \begin{array}{ccc}
                      0 &   0      &  0   \\

                        a^2\omega\wedge \omega^1&   0   &  0  \\

                           0  &  -\frac{1}{a^2}\omega\wedge \omega^2 &   0
               \end{array} \right ).
$$
The fact that the terms  $\tau^i$ do not vanish reflects the fact that this strict path structure does not appear in 
our theorem (see Lemma \ref{lemma:reductioncourbure}).  Indeed one can easily show that its automorphism group coincides with $\SU(2)$.

\vspace{.5cm}

{\bf Example:}  \label{model:SL(2,R)} Let $\SL(2,\R)$ with its left invariant vector fields defined by a Lie algebra basis $(E,F,H)$ with 
$[E,F]=H$ as in section \ref{subsection-SL2strict} .   The Maurer-Cartan equations, for a dual basis  $\omega^1,\omega^2,\omega$, are:

$$d \omega+\omega^1\wedge \omega^2=0,$$
$$d\omega^1 +2 \omega^1 \wedge \omega=0,$$
$$
d\omega^2 -2\omega \wedge \omega^2=0.
$$
Comparing with Proposition \ref{theta11} again we obtain that, for the strict path structure defined by $(\R E,\R F,H)$, $\tau^1 = \tau^2= 0$ and locally, the pull-back by a section of the component $w$ of the connection is $w=\frac{2}{3}\omega$.  The  curvature form of the connection is then given by
$$
\Omega=\left ( \begin{array}{ccc}
                      \frac{2}{3}\omega^2\wedge\omega^1 &   0      &  0   \\

                        0&   -\frac{4}{3}\omega^2\wedge\omega^1   &  0  \\

                           0  & 0 &   \frac{2}{3}\omega^2\wedge\omega^1
               \end{array} \right ).
$$
One can think of $\SL(2,\R)$ with the above strict path structure as  a constant curvature  model. Observe that one can vary the curvature by choosing  different multiples of $H$.  The curvature sign corresponds then to different choices of orientation.

\vspace{.5cm}

Bianchi identities are obtained differentiating the structure equations. They are described in the following equations:

 \begin{equation}\label{dtheta11}
d w=R\omega^1\wedge \omega^2+W^1\omega\wedge\omega^1+W^2\omega\wedge\omega^2
\end{equation}
\begin{equation}\label{dtau1}
d\tau^1+3\tau^1\wedge w=3W^2\omega^1\wedge \omega^2+S^1_1\omega\wedge\omega^1+S^1_2\omega\wedge\omega^2
\end{equation}
\begin{equation}\label{dtau2}
d\tau^2-3\tau^2\wedge w=3W^1\omega^1\wedge \omega^2+S^2_1\omega\wedge\omega^1+S^2_2\omega\wedge\omega^2
\end{equation}
Moreover, we have the relation
$$
S^1_1=S^2_2=\tau^1_2\tau^2_1.
$$
\begin{rmk}\label{rmkbianchi}
If $\tau^1=0$ (respectively $\tau^2=0$) on an open set, 
then the Bianchi identities imply that $W^2=0$ ($W^1=0$) on this open set.
More precisely, equations \eqref{dtau1} and \eqref{dtau2} shows that $W^1$ (respectively $W^2$)
is determined by $\tau^2$ (resp. $\tau^1$).
\end{rmk}

\section{Proof of Theorem \ref{thm:GhyssansAnosov}}\label{subsection:preuvetheoremeprincipal}
From now on and in all this section, 
we consider a strict path structure $\Tm=(E^1,E^2,\theta)$
on a compact three-dimensional manifold $M$,
and we denote by $\mathcal{C}=(\hat{M},\varpi)$ its normal Cartan geometry constructed in the previous paragraph,
where $\pi\colon\hat{M}\to M$ is a $P$-principal bundle over $M$ and $\varpi\colon\hat{M}\to\mathfrak{g}$
a $\mathfrak{g}$-valued Cartan connection on $\hat{M}$.
We also denote by $K\colon\hat{M}\to\mathcal{W}$ the curvature map of $\mathcal{C}$ (see \eqref{equationW}).

\subsection{Curvature types}\label{soussection:typecourbure}
With the notations of \eqref{equationW},
we introduce the $P$-invariant line 
\[
\mathcal{D}\coloneqq
\enstq{K\in\mathcal{W}}{
W^1=\tau^2_1=W^2=\tau^1_2=0
}
\]
of the space $\mathcal{W}$, and we say that a strict path structure $(M,\Lm)$ is of \emph{type $\mathcal{D}$}
if its curvature has values in $\mathcal{D}$.
\begin{lem}\label{lemma:reductioncourbure}
Let $(M,\Tm)$ be a compact connected strict path structure of class $\cc^2$
having a non-compact automorphism group.
\begin{enumerate}
\item $(M,\Tm)$ is of type $\mathcal{D}$.
\item Moreover if $\Tm$ is of class $\cc^3$ or
has a dense $\Aut^{loc}$-orbit, then its curvature map is constant.
\end{enumerate}
\end{lem}
\Pf
We recall that the automorphism group $\Aut(\Tm)$ of $\Tm$ preserves a volume form $\mu_{\Tm}=\theta\wedge d\theta$ on $M$.
We will say that a point $x\in M$ is \emph{recurrent} for the action of $\Aut(\Tm)$,
if there exists a sequence $(f_k)$ in $\Aut(\Tm)$, going to infinity in the sense that $(f_k)$ eventually
leaves every compact set (of $\Aut(\Tm)$ for the compact-open topology) 
once and for all (which will be denoted by $f_k\to\infty$), 
and such that $(f_k(x))$ converges to $x$. 
As $M$ is compact, $\mu_\Tm(M)$ is finite and Poincarr\'e recurrence Theorem applies:
since $\Aut(\Tm)$ is non-compact,
$\mu_\Tm$-almost every point of $M$ is recurrent for the action of $\Aut(\Tm)$
(see for instance \cite[Theorem 2.2.6]{FK}).
In particular, the set of recurrent points is dense in $M$. \\
1) Let $x\in M$ be a recurrent point for the action of $\Aut(\Tm)$, and $f_k\to\infty$ in $\Aut(\Tm)$ such that $f_k(x)\to x$.
Let us assume by contradiction that $\Tm$ is not of type $\mathcal{D}$ at $x$, \emph{i.e.} that
for some $\hat{x}_0\in\pi^{-1}(x)$, one of the components of $K(\hat{x}_0)$ different from $R$ is non-zero.
To fix the ideas, we assume that $\tau^1_2(\hat{x}_0)\neq0$  (the proof being the same in the other cases, \emph{mutatis mutandis}).
Introducing
\[
\hat{M}_0\coloneqq
\enstq{\hat{x}\in\hat{M}}{\tau^1_2(\hat{x})=\tau^1_2(\hat{x}_0)},
\]
the action \eqref{equationactionAWT} of $P$ on $\mathcal{W}$ shows that
for any $\hat{x}\in\hat{M}$, if $\tau^1_2(\hat{x})\neq0$ then the fiber of $\hat{x}$ meets $\hat{M}_0$.
Consequently, $M_0\coloneqq\pi(\hat{M}_0)$ is open and $\hat{M}_0$ is a reduction of the restriction of $\pi\colon\hat{M}\to M$ to $M_0$,
whose structural group is finite according to \eqref{equationactionAWT} 
(of cardinal $1$ or $2$ depending on which component is non-zero).
Since $(f_k(x))$ converges to $x\in M_0$, there exists $\hat{x}_k\in\hat{M}_0\cap\pi^{-1}(f_k(x))$ converging to $\hat{x}_0$.
But the $f_k$ being automorphisms of $\mathcal{C}$, their lifts $\hat{f}_k$ in $\hat{M}$ preserve $K$, and thus $\hat{M}_0$.
Hence for any $k$, $\hat{f}_k(\hat{x}_0)$ and $\hat{x}_k$ are in the same fiber of the finite principal bundle $\hat{M}_0$.
Since $(\hat{x}_k)$ converges in $\hat{M}_0$, this forces some subsequence of $(\hat{f}_k(\hat{x}))$ to do the same.
But $M$ is connected and the $\hat{f}_k$ preserve the parallelism $\varpi$, so the corresponding subsequence of $(f_k)$ is convergent,
which contradicts $f_k\to\infty$.
Finally, $K$ has values in $\mathcal{D}$ above any recurrent point for $\Aut(\Tm)$,
and thus everywhere by density of the recurrent points in $M$, and continuity of $K$. \\
2) The curvature of $\Tm$ is thus reduced to the continuous real-valued function $R$,
which is constant on each fiber of $\hat{M}$ according to the
equation \eqref{equationactionAWT},
and preserved by local automorphisms of $(M,\Tm)$.
If $\Tm$ has a dense $\Aut^{loc}$-orbit, then the claim immediately follows.
We now assume that $\Tm$ is of class $\cc^3$,
and we still denote by $R\colon M\to \R$ the map induced
by the only remaining curvature component on $M$,
which is of class $\cc^1$.
By the density of recurrent points for $\Aut(M,\Tm)$ \cite[Theorem 2.2.6]{FK},
it is thus sufficient to show that the differential of $R$ vanishes
at recurrent points to ensure that $R$ is constant,
finishing the proof of the Lemma.
Let $x\in M$ be recurrent, \emph{i.e.} there exists
a sequence $f_k$ of automorphisms such that $f_k\to\infty$ and
$\lim f_k(x)=x$.
Let $s=(\theta,\theta^1,\theta^2)\colon U\to\hat{M}$
be a local section of $\pi$ around $x$,
and let denote
$dR=R_0\theta+R_1\theta^1+R_2\theta^2$
in this local coframe field.
According to Lemma \ref{lemma:ReebKilling}, we already know that
$R$ is a Killing vector field and thus that $R_0\equiv0$.
There exists $a_k\in\R^*$ such that
\begin{equation}\label{equation:holonomie}
(f_k^*(\theta,\theta^1,\theta^2))_x=(\theta,a_k\theta^1,a_k^{-1}\theta^2)_x,
\end{equation}
which is equivalent to
$\hat{f}_k(s(x))=s(f_k(x))\cdot a_k$.
Since $R\circ f_k=R$,
\eqref{equation:holonomie} implies
$(dR)_x=d(R\circ f_k)_x=(f_k^*dR)_x=
R_0(f_k(x))\theta_x+R_1(f_k(x))a_k\theta^1_x+R_2(f_k(x))a_k^{-1}\theta^2_x$,
and thus
$R_1(f_k(x))=R_1(x)a_k^{-1}$ and $R_2(f_k(x))=R_2(x)a_k$.
Now if $R_1(x)\neq0$ or $R_2(x)\neq0$,
this implies $\lim a_k=1$ since $\lim f_k(x)=x$,
hence that the sequence $f_k$ is convergent which
contradicts the assumption.
Finally $R_1(x)=R_2(x)=0$.
\EPf

\begin{rmk}
Note that according to the remark \ref{rmkbianchi}, it actually suffices to prove that $\tau_1^2=\tau_2^1=0$ to obtain this Lemma.
\end{rmk}

\subsection{Constant curvature models and mutations}\label{section:mutations}

The goal of this section is to describe the local geometry of $\Tm$, by proving the following result.
\begin{prop}\label{prop:GXstructure}
Let $(M,\Tm)$ be a compact connected strict path structure having a non-compact automorphism group 
and a dense $\Aut^{loc}$-orbit.
Then: 
\begin{enumerate}
\item either $\Tm$ is, up to a constant multiplication of its contact form,
induced by a $(\SL(2,\R)\times A,\SL(2,\R))$-structure on $M$, 
\item or $\Tm$ is induced by a $(\mathbf{Heis}(3)\rtimes P,\mathbf{Heis}(3))$-structure on $M$.
\end{enumerate}
\end{prop}
In this statement, we denote by $A$ the subgroup of diagonal matrices of $\SL(2,\R)$
and we use the definitions of paragraphs \ref{subsection-flatstrict} and \ref{subsection-SL2strict} 
for the respective actions of  
$\SL(2,\R)\times A$ on $\SL(2,\R)$ and $\mathbf{Heis}(3)\rtimes P$ on $\mathbf{Heis}(3)$. 

\begin{rmk}
One should note that Lemma \ref{lemma:reductioncourbure} immediately implies that the only local models are the constant curvature ones in Example \ref{model:SL(2,R)},  together with the flat one.  Indeed, if  the curvatures are constant the structures are determined.  But, in the following, we also give an exposition of the constant curvature models using mutations (see \cite{Sh} Section 6).
\end{rmk} 

\par Let us briefly recall the notion of $(G,X)$-structure. 
Let $G$ be a group of diffeomorphisms of a manifold $X$, verifying the following condition: 
for any $g$, $h$ in $G$, if the action of $g$ and $h$ coincide on some non-empty open subset of $X$, then $g=h$
(this is clearly verified for the models of Proposition \ref{prop:GXstructure}, the action of $G$ on $X$ being real analytic).
A $(G,X)$-structure on $M$ is the data of a maximal atlas of charts from $M$ to $X$ whose transition maps are restrictions
of left translations by elements of $G$.
If $G$ preserves on $X$ a strict path structure $\Tm_X$,
then any $(G,X)$-structure $\mathcal{S}$ on a manifold $M$ \emph{induces} on $M$ a strict path structure $\mathcal{T}_{\mathcal{S}}$,
by saying that the $(G,X)$-charts of $\mathcal{S}$ are local isomorphisms from $\Tm_{\mathcal{S}}$ to $\Tm_X$.
In the proposition \ref{prop:GXstructure}, $\mathbf{Heis}(3)$ (respectively $\SL(2,\R)$) is implicitly endowed with its 
$\mathbf{Heis}(3)\rtimes P$-invariant (resp. $\SL(2,\R)\times A$-invariant) strict path structure,
described in paragraph \ref{subsection-flatstrict} (resp. \ref{subsection-SL2strict}).
\par Now, we assume that $X$ is the model of a Cartan geometry, 
\emph{i.e.} that $X=G/H$ with $H$ a closed subgroup of a Lie group $G$, such that $X$ is connected
and the action of $G$ on $X$ is faithful.
This is satisfied for the models of Proposition \ref{prop:GXstructure}, with $G/H$ being respectively $(\mathbf{Heis}(3)\rtimes P)/P$,
or $(\SL(2,\R)\times A)/\Delta$ with $\Delta=\enstq{(a,a)}{a\in A}$.
Let us denote by $\mathcal{C}_X$ 
the Cartan geometry defined on $X$ by the canonical projection $G\to G/H=X$ and the Maurer-cartan connection on $G$.
Any $(G,X)$-structure $\mathcal{S}$ on a manifold $M$ \emph{induces}
a flat Cartan geometry $\mathcal{C}_{\mathcal{S}}$ on $M$ (unique up to isomorphism), 
by saying that its $(G,X)$-charts are local isomorphisms of Cartan geometries 
from $(M,\mathcal{C}_{\mathcal{S}})$ to $(X,\mathcal{C}_X)$.\footnote{The transition maps of the maximal
$(G,X)$-atlas defining the transition maps of the Cartan bundle.}
One of the principal interests of Cartan geometries lies in the converse of this statement:
any flat Cartan geometry modelled on $G/H$ is \emph{induced} by a (unique) $(G,X)$-structure.
This claim is proved for instance in \cite[Chapter 5 Theorems 5.1 and 5.2]{Sh}. \\

\par From now on and until the end of paragraph \ref{subsection:preuvetheoremeprincipal}, we prove the Proposition \ref{prop:GXstructure}.
According to Lemma \ref{lemma:reductioncourbure}, the curvature of the normal Cartan geometry $\mathcal{C}=(\hat{M},\varpi)$ of
$\Tm$ is constant equal to some $K\in\mathcal{D}$, 
whose coordinate is denoted by $R\in\R$ with the notations \eqref{equationW}.
If $R=0$, $\mathcal{C}$ is flat and
we saw previously that $\Tm$ is then induced by a $(\mathbf{Heis}(3)\rtimes P,\mathbf{Heis}(3))$-structure on $M$.
We assume from now on that $R\neq 0$, and we prove that $\Tm$ 
is induced by a $(\SL(2,\R)\times A,\SL(2,\R))$ structure on $M$.
\par The Lie algebra of $G=\mathbf{Heis}(3)\rtimes P$ is the semi-direct product $\mathfrak{g}=\mathfrak{heis}(3)\rtimes\mathfrak{p}$
where $\mathfrak{p}$ is the Lie algebra of $P$,
whose adjoint action on $\mathfrak{heis}(3)$ is described by
\begin{equation}\label{equation:bracketg}
[D,X]=X,[D,Y]=-Y,[D,Z]=0,
\end{equation} 
with $D$ the generator of $\mathfrak{p}$ verifying
\[
\exp(tD)=
\begin{pmatrix}
e^{\frac{t}{3}} & 0 & 0 \\
0 & e^{-2\frac{t}{3}} & 0 \\
0 & 0 & e^{\frac{t}{3}}
\end{pmatrix}.
\]
\par We define a bilinear $\mathfrak{g}$-valued application $[\cdot,\cdot]'$ on $\mathfrak{g}$ 
by the following relation for $u,v\in\mathfrak{g}\times\mathfrak{g}$:
\begin{equation}\label{equation:definitionmutation}
[u,v]'\coloneqq [u,v]-K(u,v).
\end{equation}
While it is clear that $[\cdot,\cdot]'$ is skew-symmetric, 
it is proved in \cite[Chapter 5 Proposition 6.8]{Sh}
(using the Bianchi identity verified by the Cartan curvature)
that $[\cdot,\cdot]'$ actually verifies the Jacobi identity, \emph{i.e.} is a new Lie bracket on $\mathfrak{g}$.
We denote by $\mathfrak{g}'$ the Lie algebra defined by the vector space $\mathfrak{g}$ endowed with the Lie bracket $[\cdot,\cdot]'$.
Denoting $r\coloneqq \frac{3R}{2}\neq0$,
$\mathfrak{g}'$ is described by the following relations:
\begin{equation}\label{equation:crochetsmutation}
[X,Y]'=Z-2rD,[D,X]'=X,[D,Y]'=-Y,[Z,X]'=[Z,Y]'=[D,Z]'=0.
\end{equation}
\par The Lie algebra of $\SL(2,\R)\times A$ is the direct sum $\mathfrak{sl}_2\oplus\mathfrak{a}$,
and the copy of
$
H=
\left(
\begin{smallmatrix}
1 & 0 \\
0 & -1
\end{smallmatrix}
\right)
$
in the right factor of $\mathfrak{sl}_2\oplus\mathfrak{a}$ is denoted by $T$.
Note that $[T,\cdot]=0$ on $\mathfrak{sl}_2\oplus\mathfrak{a}$.
We define a group isomorphism $\Lambda\colon P\to \Delta$ by
\begin{equation}\label{equation:Lambda}
\Lambda\left(
\begin{pmatrix}
a^{\frac{1}{3}} & 0 & 0 \\
0 & a^{-\frac{2}{3}} & 0 \\
0 & 0 & a^{\frac{1}{3}}
\end{pmatrix}
\right)=
\left(
\begin{pmatrix}
a^{\frac{1}{2}} & 0 \\
0 & a^{-\frac{1}{2}} \\
\end{pmatrix}
\begin{pmatrix}
a^{\frac{1}{2}} & 0 \\
0 & a^{-\frac{1}{2}} \\
\end{pmatrix}
\right).
\end{equation}
Denoting by $\varepsilon\in\{\pm1\}$ the sign of $r$, we define 
a vector space isomorphism $\lambda\colon\mathfrak{heis}(3)\rtimes\mathfrak{p}\to\mathfrak{sl}_2\oplus\mathfrak{a}$ by:
\begin{equation}\label{equation:lambda}
\lambda(X)=\sqrt{\abs{r}}E,\lambda(Y)=-\varepsilon\sqrt{\abs{r}}F,\lambda(Z)=rT,\lambda(D)=\frac{1}{2}(H+T).
\end{equation}
\begin{lem}\label{lemma:mutationgl2}
\begin{enumerate}
\item $\lambda$ is a Lie algebra isomorphism from $\mathfrak{g}'$ to $\mathfrak{sl}_2\oplus\mathfrak{a}$.
\item The differential of $\Lambda$ at the identity coincide with $\lambda\restreinta_{\mathfrak{p}}\colon\mathfrak{p}\to\mathfrak{a}$.
\item For any $u,v\in\mathfrak{g}$: $[\lambda(u),\lambda(v)]=\lambda([u,v])$ modulo $\mathfrak{a}$
(the brackets being respectively in $\mathfrak{sl}_2\oplus\mathfrak{a}$ and $\mathfrak{heis}(3)\rtimes\mathfrak{p}$).
\item For any $p\in P$: $\lambda\circ Ad_p=Ad_{\Lambda(p))}\circ\lambda$
(the adjoint actions being respectively within the Lie groups $\mathbf{Heis}(3)\rtimes P$ and $\SL(2,\R)\times A$).
\end{enumerate}
\end{lem}
\Pf
1. This a straightforward verification from the Lie brackets relations \eqref{equation:crochetsmutation}. \\
2. This directly follows from the definitions of $\Lambda$ and $\lambda$. \\
3. In fact for $u,v\in\mathfrak{g}$, $K(u,v)\in\mathfrak{p}$ and $[\lambda(u),\lambda(v)]=\lambda([u,v]')=\lambda([u,v])-\lambda(K(u,v))$,
the first equality using the first claim of the Lemma. \\
4. For $p\in P$ written as in \eqref{equation:Lambda}, the matrix of
$Ad_p$ in the basis $(X,Y,Z,D)$ is the diagonal matrix $[a,a^{-1},1,0]$, and the matrix of $Ad_{\Lambda(p)}$ in the basis $(E,F,H,T)$ is
the diagonal matrix $[a,a^{-1},1,1]$. 
The claim directly follows from the definition of $\lambda$.
\EPf

We emphasize that in order to define a Cartan geometry modeled on $G/P$, the global model $G/P$ is in fact not necessary:
the triplet $(\mathfrak{g},\mathfrak{p},P)$ along with a morphism $Ad\colon P\to\Aut(\mathfrak{g})$ 
extending the adjoint representation $Ad_{\mathfrak{p}}\colon P\to\Aut(\mathfrak{p})$ of $P$,
are in fact sufficient to define Cartan geometries and study all their properties
(see \cite[Chapter 5 Definitions 1.1 and 3.1]{Sh}).
Denoting by $\mathfrak{d}$ the Lie algebra of $\Delta$,
we consider on $M$ the Cartan geometry $\mathcal{C}'$ modeled on $(\mathfrak{sl}_2\oplus\mathfrak{a},\mathfrak{d},\Delta)$,
whose $\Delta$-principal bundle is simply $\hat{M}$ endowed with the action of $\Delta$ given by $\Lambda^{-1}$,
and whose Cartan connection is $\varpi'\coloneqq\lambda\circ\varpi$.
It is not difficult to verify, from the properties of Lemma \ref{lemma:mutationgl2}, that $\mathcal{C}'$ is indeed a Cartan geometry.
The reader will find details in \cite[Chapter 5 Proposition 6.3]{Sh}, where Sharpe calls a map $\lambda$ verifying the properties of
Lemma \ref{lemma:mutationgl2}, a \emph{model mutation} from $(\mathfrak{heis}(3)\rtimes\mathfrak{p},\mathfrak{p},P)$ to 
$(\mathfrak{sl}_2\oplus\mathfrak{a},\mathfrak{d},\Delta)$.
\par Moreover, the Lie bracket $[\cdot,\cdot]'$ was precisely defined so that the curvature map
$K'\colon\hat{M}\to\Lin(\Lambda^2(\mathfrak{sl}_2\oplus\mathfrak{a}),\mathfrak{sl}_2\oplus\mathfrak{a})$ 
of $\mathcal{C}'$ vanishes.
In fact, for $u',v'\in\mathfrak{sl}_2\oplus\mathfrak{a}$, let introduce the $\varpi'$-constant vector fields 
$U$ and $V$ on $\hat{M}$, defined by $\varpi'(U)\equiv u'$ and $\varpi'(V)\equiv v'$.
Note that, with $u=\lambda^{-1}(u')$ and $v=\lambda^{-1}(v')$, $U$ and $V$ are also the $\varpi$-constant vector fields associated
to $u'$ and $v'$: $\varpi(U)\equiv u$ and $\varpi(V)\equiv v$.
The definition of $K'$ and the Cartan formula for $d\varpi'$ give
$K'(u',v')=[u',v']-\varpi'([U,V])$, and also $K(u,v)=[u,v]_{\mathfrak{g}}-\varpi([U,V])$.
Therefore, since $\lambda$ is a Lie algebra isomorphism from $\mathfrak{g}'$ to $\mathfrak{sl}_2\oplus\mathfrak{a}$,
the definition of $[\cdot,\cdot]'$ gives
$K'(u,v)=\lambda([u,v]_{\mathfrak{g}'})-\lambda\circ\varpi([U,V])=\lambda([u,v]_{\mathfrak{g}}-K(u,v)-\varpi([U,V]))=0$.
\par The Cartan geometry $\mathcal{C}'$ being flat, it is induced by a $(\SL(2,\R)\times A,\SL(2,\R))$-structure on $M$
(see for instance \cite[Chapter 5 Theorems 5.1 and 5.2]{Sh}).
This $(\SL(2,\R)\times A,\SL(2,\R))$-structure induces on $M$ a path structure 
$\mathcal{U}=(F^1,F^2,\theta')$ whose Reeb vector field is denoted by $R'$, which is characterized by
\[
\Diff{}{\pi}\circ\varpi'^{-1}(\R E)=F^1, \Diff{}{\pi}\circ\varpi'^{-1}(\R F)=F^2, \Diff{}{\pi}\circ\varpi'^{-1}(H)=R'.
\]
Since $\varpi'=\lambda\circ\varpi$, and $\lambda^{-1}(H)=2D-\frac{1}{r}Z$,
we thus have
$\Diff{}{\pi}\circ\varpi^{-1}(\R X)=F^1$, $\Diff{}{\pi}\circ\varpi^{-1}(\R Y)=F^2$, 
$\Diff{}{\pi}\circ\varpi^{-1}(Z)=-rR'$ (because $\Diff{}{\pi}\circ\varpi^{-1}(D)=0$).
According to the link \eqref{equationrelationCT} between the normal Cartan geometry $\mathcal{C}$ 
and its induced strict path structure $\Tm=(E^1,E^2,\theta)$, 
$(F^1,F^2,\theta')$ is thus equal to $(E^1,E^2,-r\theta)$.
In other words, the $(\SL(2,\R)\times A,\SL(2,\R))$-structure of $M$ 
is associated with the original strict path structure $\Tm$, up to a constant multiplication of its contact form,
which concludes the proof of Proposition \ref{prop:GXstructure}.

\subsection{Completeness of the structure}\label{soussection:completude}

\par A classical fact about $(G,X)$-structures is the existence of a \emph{developping map} $\delta\colon\tilde{M}\to\tilde{X}$
between the respective universal covers of $M$ and $X$, which is a local diffeomorphism extending the $(G,X)$-charts of the atlas.
In our case, $\delta$ is a local isomorphism of strict path structures, and
$\tilde{X}$ is either $\mathbf{Heis}(3)$ or $\widetilde{\SL}(2,\R)$.
Denoting by $\tilde{G}$ the group 
$\mathbf{Heis}(3)\rtimes P$ in the first case and $\widetilde{\SL}(2,\R)\times\tilde{A}$ in the second one,
$\delta$ is equivariant with respect to a morphism $\rho\colon\pi_1(M)\to\tilde{G}$ called the holonomy morphism.
If $\delta$ is a covering map, then the $(G,X)$-structure is said to be \emph{complete}
and $\delta$ is a diffeomorphism. In this case $\Gamma\coloneqq\rho(\pi_1(M))$ acts freely, properly and 
cocompactly on $\tilde{X}$, yielding an isomorphism of $(G,X)$-structures (and thus of strict path structures) 
between $M$ and $\Gamma\backslash\tilde{X}$
(for more details about $(G,X)$-structures, see for instance \cite{T}).
\par To conclude the proof of Theorem \ref{thm:GhyssansAnosov},
we are thus left to understand that the $(G,X)$-structure of the compact manifold $M$
is complete. 
This completeness follows from previous works about Lorentzian metrics in \cite{DZ} and \cite{Kl}, that we now present.
{We then give an independent argument, in a specific dynamical situation.} \\

\noindent{\em End of the proof of Theorem \ref{thm:GhyssansAnosov}}
If $G$ preserves on $X$ a Lorentzian metric $g_0$, the $(G,X)$ structure of $M$ induces a metric $g$ for which $\delta$ is locally isometric.
If $g$ is complete, then $\delta$ is a covering map
(and if $g_0$ is complete, the converse is true).
\par The group $\mathbf{Heis}(3)\rtimes P$ preserves on $\mathbf{Heis}(3)$ the (non-flat) left-invariant Lorentzian metric $g_0$
called \emph{Lorentz-Heisenberg} and described as follows on its Lie algebra $\mathfrak{heis}(3)$:
$Z$ has norm $1$, $\Vect(X,Y)$ is orthogonal to $X$, $X$ and $Y$ are isotropic vectors, and $g_0(X,Y)=1$.
It is shown in \cite[\S 8.1]{DZ} that any compact Lorentzian manifold locally isometric to the Lorentz-Heisenberg metric is complete, 
and that if a quotient $\Gamma\backslash\mathbf{Heis}(3)$ is compact
with $\Gamma\subset\mathbf{Heis}(3)\rtimes P$, then $\Lambda\coloneqq\Gamma\cap\mathbf{Heis}(3)$ has finite index in $\Gamma$.
In other words, any compact $(\mathbf{Heis}(3)\rtimes P,\mathbf{Heis}(3))$-structure is, up to a finite covering of $M$, a quotient
$\Gamma\backslash\mathbf{Heis}(3)$ with $\Gamma$ a cocompact lattice of $\mathbf{Heis}(3)$, 
which concludes the proof in the case of $\mathbf{Heis}(3)$.
\par The left-invariant Lorentzian metric of $\widetilde{\SL}(2,\R)$ defined by the Killing form on $\mathfrak{sl}(2)$
is the universal cover of the \emph{Anti-de-Sitter} space (local model of curvature $-1$ of Lorentzian metrics)
and is preserved by the action of $\widetilde{\SL}(2,\R)\times\tilde{A}$.
Klingler showed in \cite{Kl} that any compact Lorentzian manifold locally isometric to Anti-de-Sitter is complete,
which concludes the proof in the case of $\widetilde{\SL}(2,\R)$, and hence the proof of Theorem \ref{thm:GhyssansAnosov}. 
\EPf

\par Let us give an alternative and easier argument for the completeness of the structure, independent from \cite{DZ} and \cite{Kl},
under the following stronger dynamical hypothesis.
\begin{center}
$(H)$ For any $x\in M$ and $i\in\{1,2\}$,
there exists a sequence $(f_k)$ of automorphisms 
of $\Tm=(E^1,E^2,\theta)$ such that $\norme{\Diff{x}{f_k}\restreinta_{E^i}}$ converges to $0$
(for some Riemannian metric on $M$).
\end{center}
We saw previously that the $(G,X)$-structure of $M$ induces a Lorentzian metric invariant by the automorphism group of $\Tm$,
and we denote by $\nabla$ its Levi-Civita connection.
Note that any geodesic tangent to one the three directions $E^1$, $E^2$ or $R$ (the Reeb vector field of $\theta$)
stays tangent to this direction (because this is the case in the model $X$).
In particular, $R$ being complete since $M$ is compact, this shows that the geodesics of $\nabla$ tangent to $R$ are complete.
The following Lemma already appears in \cite{BFL}.
\begin{lem}\label{lemma:completudealphabeta}
The geodesics of $\nabla$ tangent to the direction $E^1$ and to the direction $E^2$ are complete.
\end{lem}
\Pf
We endow $M$ with a Riemannian metric.
By compacity of $M$, there exists a $\varepsilon>0$ such that for any $x\in M$ and $u$ a unit vector in $E^1(x)$ or $E^2(x)$,
the geodesic starting from $x$ with speed $u$ is defined until time $\varepsilon$.
For $x\in M$, $i\in\{1,2\}$ and $u\in E^i(x)\setminus\{0\}$, there exists by hypothesis  a sequence $(f_k)$ of automorphisms 
of $\Tm$ such that $\norme{\Diff{x}{f_k}(u)}$ converges to $0$,
and there exists thus $k$ such that $\norme{\Diff{x}{f_k}(u)}<\varepsilon$.
The geodesic at $f_k(x)$ with speed $\Diff{x}{f_k}(u)$ is thus defined until time $1$, 
and its image by $f_k^{-1}$ also, because $f_k$ is an affine application.
But this is the geodesic starting from $x$ with speed $u$, which concludes the proof.
\EPf
\par Thanks to this result, we now prove the completeness of the $(G,X)$-structure under the hypothesis $(H)$.
We will say that a path $\gamma\colon\intervalleff{0}{1}\to X$ starting at $x=\gamma(0)$ \emph{lifts} to $\tilde{M}$ 
\emph{from} $\tilde{x}\in\delta^{-1}(x)$, if there exists $\tilde{\gamma}\colon\intervalleff{0}{1}\to\tilde{M}$
(called the lift of $\gamma$ from $\tilde{x}$)
such that $\tilde{\gamma}(0)=\tilde{x}$ and $\delta\circ\tilde{\gamma}=\gamma$.
Let us denote by $\Fitan{\tilde{X}}=E^1_0\oplus E^2_0 \oplus \R R_0$ the smooth splitting defined by the strict path structure of $\tilde{X}$,
$R^0$ being its Reeb vector field. 
For $x\in\delta(\tilde{M})$, $\tilde{x}\in\delta^{-1}(x)$ and $i\in\{1,2\}$, 
since the geodesic tangent to $E^i$ starting from $\tilde{x}$ is complete,
any path starting from $x$ and tangent to $E^i_0$ lifts from $\tilde{x}$.
According to \cite[Lemma 7.2]{Mm}, 
this implies that $\delta$ satisfies the path-lifting property, and is thus a covering from $\tilde{M}$ to $\tilde{X}$.

\appendix
\section{Proof of Proposition \ref{po:pasHeis}}\label{subsection:pasHeisenberg}
Let $\Gamma\subset G=\mathbf{Heis}(3)\rtimes P$ be a discrete subgroup, acting freely, properly and cocompactly on $\mathbf{Heis}(3)$
by the action described in paragraph \ref{subsection-flatstrict},
and let us assume by contradiction that the central flow acting on $\Gamma\backslash\mathbf{Heis}(3)$ is not periodic.
We will denote by $\mathcal{Z}=\{z^t\}_{t\in\R}$ the one-parameter center of $\mathbf{Heis}(3)$.
\par Our hypothesis implies that $\Gamma\cap\mathcal{Z}=\{e\}$,
otherwise if $z^{t_0}\in\Gamma\cap\mathcal{Z}$ with $t_0\neq0$, then for any $n\in\N$ and $x\in\mathbf{Heis}(3)$
we would have $z^{nt_0} \Gamma x=\Gamma x$, \emph{i.e.} the central flow acts periodically.
On the other hand, $\Lambda\coloneqq\Gamma\cap\mathbf{Heis}(3)$ 
is the kernel of the restriction to $\Gamma$ of the second projection $\pi_2\colon hp\in \mathbf{Heis}(3)\rtimes P\mapsto p\in P$.
Let us assume by contradiction that $\Lambda=\{e\}$. 
Then $\Gamma$ is isomorphic to $\pi_2(\Gamma)\subset A$ and is thus abelian.
Since algebraic subgroups have a finite number of connected components, 
the identity component $(\overline{\Gamma}^Z)^0$ of the Zariski closure of $\Gamma$ in $G$ 
has finite index in $\overline{\Gamma}^Z$,
and $\Gamma'\coloneqq\Gamma\cap(\overline{\Gamma}^Z)^0$ is thus a finite index subgroup of $\Gamma$. 
But $(\overline{\Gamma}^Z)^0$ is a (non-trivial) closed connected abelian subgroup of $G$,
and as such, one easily verifies that it is isomorphic to $\R$ or $\R^2$.
Indeed, abelian Lie subalgebras of $\mathfrak{g}$ of dimension strictly greater than one,
are either contained in $\mathfrak{heis}(3)$ and therefore isomorphic to $\R^2$,
or conjugated to $\mathfrak{a}_1\oplus\mathfrak{z}(\mathfrak{heis}(3))$ by the adjoint action of $G$.
Consequently, $\Gamma'$ is isomorphic to $\Z$ or $\Z^2$.
But $\Gamma'$ having finite index in $\Gamma$, its action on $\mathbf{Heis}(3)$ remains proper and cocompact,
implying that $\Gamma'$ has cohomological dimension equal to $3$ since $\mathbf{Heis}(3)$ is contractible (see \cite[Prop. 8.1 p.210]{Br}).
This contradicts the fact that $\Gamma'$ is isomorphic to $\Z$ or $\Z^2$ 
(that have respective cohomological dimensions equal to $1$ and $2$), and finally shows that $\Lambda$ is a 
non-trivial discrete subgroup of $\mathbf{Heis}(3)$ whose intersection with $\mathcal{Z}$ is trivial.
\par Since the commutator subgroup of $\mathbf{Heis}(3)$ is contained in its center,
$[\Lambda,\Lambda]\subset\Gamma\cap\mathcal{Z}=\{e\}$, \emph{i.e.} $\Lambda$ is abelian.
We introduce as before $\Lambda'\coloneqq\Lambda\cap(\overline{\Lambda}^Z)^0$, 
and we emphasize that this is a finite index subgroup of $\Lambda$,
as well as a normal subgroup of $\Gamma$ (indeed $\Lambda$ is normal in $\Gamma$,
and so are its Zariski closure $\Lambda$ and its identity component $(\overline{\Lambda}^Z)^0$).
Since $(\overline{\Lambda}^Z)^0$ is a closed connected abelian subgroup of $\mathbf{Heis}(3)$, there are two possibilities:
\begin{enumerate}
\item either $(\overline{\Lambda}^Z)^0$ is a one-parameter subgroup, and $\Lambda'$ is isomorphic to $\Z$,
\item or the Lie algebra of $(\overline{\Lambda}^Z)^0$ is, up to conjugation in $\mathbf{Heis}(3)$, equal to $\Vect(X,Z)$ or $\Vect(Y,Z)$, 
and $\Lambda'$ is then isomorphic to $\Z$ or $\Z^2$.
\end{enumerate}
This already shows that $\Lambda\subsetneq\Gamma$, because $\Z$ and $\Z^2$ cannot act freely, properly, cocompactly on 
the contractible manifold $\mathbf{Heis}(3)$.
So if $\Lambda'$ was not isomorphic to $\Z$, we would be in the second case, and
there would be some $\gamma\in\Lambda'$ and $g\in\Gamma$ of the form
\begin{equation}\label{equation:notationsgammag}
\gamma=[x,0,z]\coloneqq
\begin{pmatrix}
1 & 0 & 0 \\
0 & 1 & 0 \\
z & x & 1 
\end{pmatrix}
\text{~and~}
g=
\begin{pmatrix}
a & 0 & 0 \\
y' & a^{-2} & 0 \\
z' & x' & a
\end{pmatrix},
\end{equation}
with $z\neq0$ and $a\neq1$ (we assume $Lie((\overline{\Lambda}^Z)^0)=\Vect(X,Z)$, the second case being treated in the same way). 
Then we have $g^n\gamma g^{-n}=[a^{3n}x,0,z-na^{3n}xy']\in\Gamma$,
converging to $[0,0,z]$ at $+\infty$ or $-\infty$ since $a\neq1$.
But $\Gamma$ is discrete and hence closed, so $[0,0,z]\in\Gamma\cap\mathcal{Z}\setminus\{e\}$, which is a contradiction.
\par Finally, $\Lambda'$ is isomorphic to $\Z$.
With $\gamma=[x,y,z]$ a generator of $\Lambda'$ and $g\in\Gamma$ as in \eqref{equation:notationsgammag}, 
a direct calculation  gives $g\gamma g^{-1}=[a^3x,a^{-3}y,z+a^{-3}yx'-a^3xy']$.
Consequently, $g\gamma g^{-1}=\gamma$ implies $a=1$ \emph{i.e.} $g\in\mathbf{Heis}(3)$ (note that $x$ or $y$ is non-zero by hypothesis),
and the centralizer of $\Lambda'$ in $\Gamma$ is thus equal to $\Gamma\cap\mathbf{Heis}(3)=\Lambda$.
Since $\Lambda'$ is normal in $\Gamma$ and isomorphic to $\Z$, this forces $\Lambda$, and hence $\Lambda'$, 
to have finite index in $\Gamma$.
Therefore, $\Gamma'$ acts freely, properly and cocompactly on $\mathbf{Heis}(3)$, 
which contradicts the fact that $\Lambda'$ is isomorphic to $\Z$ and 
concludes the proof of Proposition \ref{po:pasHeis}.
Indeed, the centralizer of $\Lambda'$ in $\Gamma$ is the kernel of the morphism 
$\Phi\colon g\in\Gamma\mapsto\{h\mapsto ghg^{-1}\}\in\Aut(\Lambda')$.
But since $\Lambda'\simeq\Z$, $\Aut(\Lambda')=\{\pm\id\}$, and hence
$\Lambda\backslash\Gamma=(\Ker\Phi)\backslash\Gamma\simeq\Phi(\Gamma)$ has cardinal at most $2$.

\begin{flushleft}
  \textsc{E. Falbel\\
  Institut de Math\'ematiques \\
  de Jussieu-Paris Rive Gauche \\
CNRS UMR 7586 and INRIA EPI-OURAGAN \\
 Sorbonne Universit\'e, Facult\'e des Sciences \\
4, place Jussieu 75252 Paris Cedex 05, France \\}
 \verb|elisha.falbel@imj-prg.fr|
 \end{flushleft}
 
 \begin{flushleft}
  \textsc{M. Mion-Mouton\\
  Institut de Recherche Math\'ematique Avanc\'ee \\
 UMR 7501\\Universit\' e de Strasbourg\\
 7 rue Ren\' e Descartes F-67000 
 Strabourg, France\\}
 \verb|martin.mionmouton@math.unistra.fr|
 \end{flushleft}

\begin{flushleft}
  \textsc{J. M.  Veloso\\
  Faculdade de Matem\' atica - ICEN\\
Universidade Federal do Par\'a\\66059 - Bel\' em- PA - Brazil}\\
  \verb|veloso@ufpa.br|
\end{flushleft}

\end{document}